\documentclass[11pt,twoside]{amsart}

\usepackage{multirow}

\usepackage{tcolorbox}

\usepackage{graphicx}
\graphicspath{ {images/} }
\usepackage{amsmath}
\usepackage{subcaption}
\usepackage{booktabs}
\usepackage{array}
\usepackage{colortbl}
\usepackage{enumerate}
\usepackage{tabularx}
\usepackage{amssymb}
\usepackage{enumitem}
\usepackage{float}
\restylefloat{table}
\usepackage[margin=1in,a4paper]{geometry}

\tcbuselibrary{raster}
\newtheorem{theorem}{Theorem}[section]
\newtheorem{rem}[theorem]{Remark}

\title{\bf The Cognition of Counterexample in Mathematics Students}

\author{Shannon Ezzat}
\address{Department of Mathematics and Statistics\\
University of Winnipeg\\
515 Portage Avenue,\\ 
Winnipeg, Manitoba, Canada R3B 2E9}
\email{s.ezzat@uwinnipeg.ca (Shannon Ezzat)}

\author{Scott Rodney}
\address{Department of Mathematics, Physics, and Geology\\
Cape Breton University\\
Box 5300, 1250 Grand Lake Road\\
Sydney, Nova Scotia, Canada B1P 6L2}
\email{scott.rodney@gmail.com (Scott Rodney)}

\date{{}}

\begin{document}

\maketitle
\thispagestyle{empty}

\begin{abstract} Studying Mathematics requires a synthesis of skills from a multitude of academic disciplines; logical reasoning being chief among them. This paper explores mathematical logical preparedness of students entering first year university mathematics courses and also the effectiveness of using logical facility to predict successful course outcomes.  We analyze data collected from students enrolled at the University of Winnipeg in a pre-service course for high school teachers.   We do find that, being able to successfully answer logical questions, both before and after intervention, are significant in relation to improved student outcomes.

\end{abstract}

\section{Introduction}

This paper explores mathematical logic preparedness and its effect on success in the first year university mathematics classroom.  It is well known that skill sets among students entering courses in university mathematics vary widely from student to student and thus an understanding of the abilities and effects of logical aptitude is a worthwhile endeavour to help improve student success in mathematics.  

This work focuses on two themes: student facility in reasoning about counterexample, and how well assessment of this logical skill can predict successful course outcomes in introductory mathematics courses.  

To address these questions, a 3-month study of university students was initiated in MATH-2903, Math for Early/Middle Year Teachers \footnote{University of Winnipeg (UW), Manitoba, Canada} during the winter semester of the 2017-18 academic year. As indicated by its name, this class is designed for pre-service elementary- and middle-school teachers. The prerequisite for the course is any university-level math class from Grade 12 in Manitoba (with a grade above 65\% for applied math or any passing grade in precalculus) and does not count towards any math major course, or as a required math course for Bachelor of Science students. The students in this course predominantly have poor to adequate arithmetic and algebraic skills and a substantially negative attitude towards mathematics as a subject.  

Our methods are inspired by several works on student ability with logic and counterexample including \cite{BF} which discusses undergraduate students understanding of various concepts in calculus using the tool of diagnostic tests, intervention/primers, and quizzes; \cite{L}, which studies the cognition of counterexample of high school students in a similar fashion; and \cite{IS}, which studies undergraduate (and staff) success in the Wason Selection Task. We use these tools to measure students understanding of logical statements, and specifically counterexamples to these statements.

A careful description of our study is found in Section 2. Section 3 contains statistical analysis discussing evidence of predictive assessment question types. A summary of our results together with concluding remarks is found in Section 4. Note that in a future note the authors will report results for a similar study from MATH-1105, Differential and Integral Calculus I \footnote{Cape Breton University (CBU), Nova Scotia, Canada}.

\indent Reasoning logically using counterexample is linked to the idea of conditional statements.  We begin with a discussion about conditional statements and the notion of counterexample itself and its importance in mathematical pedagogy, including a particular instance of a question type designed by Peter Wason \cite{W} called the Wason Selection Task. We refer to this as the Wason question in our tests.

\subsection{Conditional Statements} A conditional statement consists of two parts, an antecedent (if-part) and a consequent (then-part). In general, these are of the form \emph{``If P then Q"}, for some statements P and Q. An example would be the following: 
   \begin{center}\emph{If an integer is even then the rightmost digit of the number is either a $2$,$4$,$6$ or $8$}. \end{center}
An important component of logic and mathematics in general, is determining the truth value (True(T) or False(F)) of a given conditional statement. To that end it is sometimes helpful to reword conditionals using a universal quantifier. In general, \emph{``If P then Q"} becomes \emph{``All P are Q"} and specific to the example above:
   \begin{center}
       \emph{All even integers have the rightmost digit being $2$,$4$,$6$ or $8$}.
   \end{center}
To falsify such a universally quantified statement is to demonstrate the existence of an object that satisfies statement P but does not satisfy statement Q; in our example, $10$ is such an object. We call such an object a \emph{counterexample} to the conditional statement. In fact, a conditional statement is false exactly when such an object exists. For more details, the reader is encouraged to refer to the introductory chapters of most textbooks in logic or discrete mathematics. It is important to see that the idea of counterexample is not necessarily intuitive. We note one more example from our diagnostic study to illustrate this point:  
   \begin{quote} If two shapes have the same perimeter then they have the same area.
   \end{quote}
On the diagnostic test,  31 of 62 responded correctly indicating that the statement is false while only 8 of those justified their conclusion using a counterexample.  Counterexamples provided by correct respondents most often used rectangles with few referring to circles or triangles.\\

\subsection{The Importance of Counterexample}
In mathematical reasoning, the notion of counterexample is fundamental. Counterexamples allow us to understand the limits of mathematical results and why each condition in a mathematical logical statement is needed to justify the conclusion. In \cite{K}, Klymchuck writes that, for undergraduate students, understanding by way of counterexamples ``can improve students’ conceptual understanding in mathematics, reduce their common misconceptions, provide a broader view on the subject and enhance students’ critical thinking skills''. For students that plan on a career in education Vinsonhaler and Lynch \cite{VL} believe that primary and secondary school teachers `` support student reasoning and thinking and promote productive struggle by incorporating counterexamples into the classroom.'' Zazkis and Chernoff \cite{ZC} suggest that counterexamples can be particularly convincing when they are ``in accord with individuals' example spaces''. \\
   
   %============
   %============

   \subsection{The Wason Question}
   
The card selection question was originally reported by Wason \cite{W}. In his experiment, there are cards on a table and the participant is told that there is a number on one side and a letter on the other. The participants are told of a rule, for example \emph{``If a card has a D on one side then it has a 3 on the other side"} or, reworded as a universally quantified statement \emph{``All cards with a D on one side have a 3 on the other side''}. Four cards are presented that correspond to each logical case (P true, P false, Q true, Q false) and the participant is asked \emph{``to select all those cards, but only those cards, which you would have to turn over in order to discover whether or not the rule has been violated."}\\
\indent To answer the question, the participant is required to assess the conditional statement four times in relation to each card presented.  The correct answer is to select the cards where either P is true or Q is false since these conditions are necessary (but not sufficient) for the card to be a counterexample or exception to the rule. Interestingly, most of the population does not make this connection.  Indeed, Inglis and Simpson in \cite{IS} found from an online, self-selected sample at the University of Warwick, only 29 percent of math students, 8 percent of history students and 49 percent of mathematics staff responded correctly.  During our study described below, only 2 students answered this style of question correctly on the pre-test; see \S 2. %Similar statistics related to simpler conditional statements were also found during our study.  

\subsection{Proving a conditional true}

More difficult perhaps is the skill of reasoning the justification of a true conditional.  Consider the following sentence:
    \begin{quote}
    If the sum of two positive whole numbers is odd then their product is even.
    \end{quote}
To justify the truth of this statement the participant must covey that, since the sum of the numbers is odd, one number must be even while the other odd.  This requires more algebraic ability, and usually a considerable amount of lecturing to convey; indeed, many universities have a course in the first or second year of the program (usually Discrete Mathematics or a standalone Introduction to Proofs course) that spend considerable time teaching the concept of correct proof writing. We note that the diagnostic test did not include any questions of this type, however as part of the course material, some lecture time was spent teaching how to prove a conditional statement for simple examples; e.g. ``If $a,b,c$ are three consecutive whole numbers then their sum has a factor of 3''.

\indent As we will discuss in more detail, the ability to understand and construct counterexamples to a conditional statement is not inherent in most students entering Math 2903. These results are not surprising in the context of similar statistical results found in \cite{IS} and \cite{L}.  However, perhaps more interesting is that on average, students who are able to respond well to these problems have better course outcomes.  We also observe that students who learn these skills through in-course intervention experience better course outcomes on average.  As such, our study suggests that logical skill assessment can be used as a tool to help improve course outcomes for students. Such assessment is also effective in identifying struggling students early on so that they can find effective help.  Lastly, it is also quite interesting that a simple intervention can on average have a measurable impact on student success.\\
   
   %=============
   %=============

\section{The Study}

The data used for this study comes from Math 2903, Math for Early/Middle Year Teachers at the University of Winnipeg.  We chose students in Math 2903 since it was, most likely, the students' first university mathematics class. Indeed, the likelihood that our student populations in the class had been exposed to logical reasoning on an academic level (for example, a course on critical thinking) was quite low. This is borne out in the results.  

The study was conducted in four phases.  In the first week of classes (Phase I), students were given a diagnostic test to benchmark their understanding of logical thinking. In the third week (Phase II), students were given an intervention (75-minute lab/lecture on the idea of conditional statements and counterexamples) complete with a worksheet for practice.  During the 6th week of classes (Phase III), a post test was given as part of the in-class midterm exam.  Final grades were tallied (Phase IV) and small sample size statistical methods were used to analyze aggregate data collected.  In what follows, we give detailed descriptions of each phase of the study but leave Phase IV to \S 3.\\

\subsection{I - The Diagnostic Test}
The Diagnostic test was designed to gauge incoming student ability to reason by counterexample and was given during the first week of class.  The problems appearing on the test were selected from 2 categories similar to those found in \cite{L}.  We describe them now.  We also note that specific questions can be found in \S 3.

\subsubsection{Type I: The Wason selection task}

As described in the introduction, this exercise concerns rule verification for a finite set of objects.  A complete description of the Wason selection task can be found in both of \cite{L} and \cite{IS}.  For this study, we depicted an array of 6 cards each with a number on one side and a letter on the other, see \eqref{card}.\\

%\definecolor{babyblue}{rgb}{0.54, 0.81, 0.94}

%\begin{figure}
\begin{eqnarray}\label{card}
\begin{tabular}{m{1cm}m{1cm}m{1cm}m{1cm}m{1cm}m{1cm}}
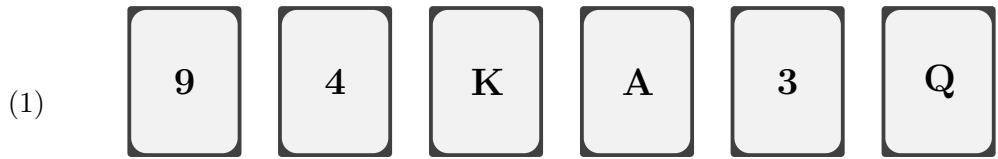
\begin{tcolorbox}[width=1.5cm,height=2cm,halign=center,valign=center,arc=2mm,outer arc=0.5mm]
\Large{\bf 9}
\end{tcolorbox}
&
\begin{tcolorbox}[width=1.5cm,height=2cm,halign=center,valign=center,arc=2mm,outer arc=0.5mm]
\Large{\bf 4}
\end{tcolorbox}
&
\begin{tcolorbox}[width=1.5cm,height=2cm,halign=center,valign=center,arc=2mm,outer arc=0.5mm]
\Large{\bf K}
\end{tcolorbox}
&
\begin{tcolorbox}[width=1.5cm,height=2cm,halign=center,valign=center,arc=2mm,outer arc=0.5mm]
\Large{\bf A}
\end{tcolorbox}
&
\begin{tcolorbox}[width=1.5cm,height=2cm,halign=center,valign=center,arc=2mm,outer arc=0.5mm]
\Large{\bf 3}
\end{tcolorbox}
&
\begin{tcolorbox}[width=1.5cm,height=2cm,halign=center,valign=center,arc=2mm,outer arc=0.5mm]
\Large{\bf Q}
\end{tcolorbox}
\end{tabular}
\end{eqnarray}
%\end{figure}

The student was told: ``Each card has a number and letter on opposite sides.  Every card with a vowel on one side has an even number on the other. Of the six cards below (it can be more than one) which must you flip over in order to determine whether or not the rule has been violated?  Try your best to explain your reasoning." Note that the correct answer to this problem is to flip the A, 9, and 3, as these could be possible counterexamples; for example, the A could have a 3 on the other side, and the 3 or 9 could have an A on the other side, and these cards would violate the rule).
%The student was asked which cards must be flipped in order to verify the rule ``If a card has a vowel on one side it must have an even number on the other".  They were also asked to justify their choices.  The correct answer to th%Interestingly, the most common response among students to this question is to check the $4$ and ace(A).  In fact, in math 1105, $XX\%$ of students answered this way.  These choices, made because $4$ is even and A is a vowel, are also common in the general population as described in a large study found in [[REF WASON]].......\\

\subsubsection{Type II: Provide a counterexample to a given conditional statement}

Here the student is asked to falsify a conditional statement by constructing a counterexample.  For example,
\begin{center}
{\quote ``If two shapes have the same perimeter then they have the same area"}
\end{center}
is seen to be false by considering two rectangles with dimensions $2\times 4$ and $5\times 1$.  This forms the required counterexample.\\

\subsection{II - The Intervention}
The intervention consisted of a 75 minute primer on conditional statements that occurred during the second week of classes.  During the intervention, students were introduced to the concept of a conditional statement, initially via a truth table, and the relationship between universally quantified statements and conditionals. Students were then shown an applied example of a liquor inspector enforcing the rule described by the conditional statement ``If you are drinking alcohol then you must be 18 years of age or older'' to reinforce the idea of when a conditional statement is false.  This was then connected to the notion of counterexample.  The intervention lecture closed with a discussion of the Wason selection task where students were shown the relationship between the truth values of the conditional sentence defining the rule and the values on each of the cards.  This was followed by a group work session where students reinforced these ideas by completing a practice worksheet with their peers. 

In addition to this 75 minute primer, as part of the course work students were given one lesson on verifying a conditional statement via an algebraic proof. The students were also given a small number of beginner practice problems, and one problem on an assignment to complete. Note that this type of question was not included on the original diagnostic test as algebraic verification is a technique that is more specialized and the instructor of the course (the first author) believed that most students would not have been able to answer the question appropriately. 

\subsection{III - The Post Test}  The post test was given as part of the MATH-2903 class midterm test. The students were given one Wason Selection Task question, as well as one question where the students would have to identify which of two conditional statements had a counterexample, provide a counterexample for that statement, and algebraically justify the other conditional statement.

\section{Data and Results}

The aim of the data analysis found in this section is to compile and compare results of the diagnostic test (DT) and post test (PT). Our results, presented in tables 1-6, analyze relationships between success/retention/difficulty with specific problems on our tests with mean final course grades.  In order to present this efficiently, we group and name specific problems of interest together with our raw data. Our data consists of evaluated student performance in terms of the following categories.\\

\begin{itemize}
\item[(i)] {\bf CR}: \emph{Correct answer and reasoning.} This evaluation is given if the student provides a clear, complete solution to the problem providing a counterexample if necessary.
\item[(ii)] {\bf C}: \emph{Correct answer with incorrect reasoning.} The student was able to answer the problem correctly but lacked sufficient justification such as, for example, an erroneous or missing counterexample.
\item[(iii)] {\bf I}: \emph{Incorrect answer.} The student has answered incorrectly with or without justification.
\item[(iv)] {\bf B}: \emph{The student left the problem blank.}
\end{itemize}

\subsection{The Diagnostic Test}
The diagnostic test consisted of 5 questions of varying type that we state now followed by student response data.  \\
\begin{itemize}
\item {\bf DT:1} If a number $n$ is less than 1, then the expression $n^2-1$ is negative.
\item {\bf DT:2} If two shapes have the same perimeter, then they have the same area.
\item {\bf DT:3} If the sum of two numbers is odd, then their product is even.
\item {\bf DT:4} If two numbers $x,y$ have positive sum, then both $x$ and $y$ are positive.
\item {\bf DT:5} This is the Wason selection task problem exactly as given in section \S 2.1.1. with cards as presented in figure \eqref{card}.
\end{itemize}

%\begin{itemize}
%\item {\bf 1(a)} If a number $n$ is less than 1, then the expression $n^2-1$ is negative.
%\item {\bf 1(b)} If two shapes have the same perimeter, then they have the same area.
%\item {\bf 1(c)} If the sum of two numbers is odd, then the product of the two numbers is even.
%\item {\bf 1(d)} If two numbers $x,y$ have positive sum, then both $x$ and $y$ are positive.
%\end{itemize}

%Question 2 was the Wason selection task mentioned above.

\begin{center}
\begin{table}[h]
\caption{Summary of student answers to counterexample questions on diagnostic test}\label{Diag}
\begin{tabularx}{0.8\textwidth}{ >{\centering\arraybackslash}X||>{\centering\arraybackslash}X 
| >{\centering\arraybackslash}X|>{\centering\arraybackslash}X 
| >{\centering\arraybackslash}X  }
     Question & CR &C &I & B  \\
     \hline
     \hline
     DT:1 & 19 & 14 & 19 & 10\\ 
     DT:2 & 8 & 23 & 20 & 11 \\
     DT:3 & 18 & 18 & 9 & 17\\
     DT:4 & 26 & 26 & 6 & 4
\end{tabularx}
\end{table}
\end{center}

For the Wason selection task problem DT:5, we analyze student responses more finely using  common observed card selections.  Category ``{A93}" refers to the student selecting the correct cards, the ace, nine, and three.  The other categories, ``A", ``A4", ``All", ``Blank" and ``Other" are self evident. \\    
\begin{center}
\begin{table}[h]\label{dtwasonresp}
\caption{Summary of student answers to Wason Selection Task question on diagnostic test}\label{Wason}
\begin{tabular}{c||c|c|c|c|c|c}
     Answers & A93 & A & A4 & All & Blank & Other  \\
     \hline
     Question 2 & 3 & 2 & 26 & 12 & 10 & 9  
\end{tabular}
\end{table}
\end{center}

\begin{rem}\label{rem1} We note that very few students correctly answered the question (A93) on the diagnostic test.  Furthermore the most common answer was to select the A and four, (A4). This choice is the common fallacy that the conditional ``if $P$ then $Q$" is equivalent to ``$P$ and $Q$".  This result aligns with previous work by Inglis and Simpson \cite{IS}, and Wason \cite{W}.
\end{rem} 

%Discussion of the other questions with references
\subsection{The Post Test:}

The post test consisted of three questions PT:1, PT:2, PT:3 on the first midterm of Math 2903. Question PT:1 was the Wason selection task problem as in \S 2.1.1, see \eqref{card}, with the altered rule 
\[\textrm{``if a card has an odd number on one side then it has a vowel on the other".}\]
The correct answer is to select the nine, three, K and Q (93KQ) as in \eqref{card}.  \\
Questions PT:2 and PT:3 appeared as two parts of a single problem that we state now. 

\medskip
\emph{One of the following conditional statements is always true and the other has a counterexample. ***Circle (a) or (b)*** of the the true statement and justify that it is true (using the direct method). Find a counterexample for the other statement (and explain why it is a counterexample). 
\begin{itemize} 
\item {\bf PT:2} (a) If $m$ and $n$ are whole numbers and $m \times n$ is even, then both $m$ and $n$ are even.
\item {\bf PT:3} (b) If $n$ is an even number then $n^2$ is a multiple of $4$. 
\end{itemize}}
\medskip

The results of the post-test questions are shown in the table below. Students are grouped into two categories: correct, and otherwise. 

\begin{center}
\begin{table}[h]
\caption{Number of students answering post-test questions correctly}\label{posttest}
\begin{tabular}{c|c|c|c}
Question & Correct Answer & Correct & Incorrect\\
\hline
PT:1 & Selected (93KQ) & 21 & 38\\
PT:2 & Provided a counterexample & 41 & 18\\
PT:3 & Provided a correct proof &12 &47\\
\hline
\end{tabular}
\end{table}
\end{center}
 
\begin{rem} In contrast with Remark \ref{rem1} and noting that this specific question type was covered in the stage II intervention, it is very encouraging that a far larger proportion of students answered correctly. However, checking the truth of ``$P$ and $Q$" is still frequent.
\end{rem}

\subsection{Relationship between counterexample knowledge and final grade}

Below are tables with mean final grades in Math 2903 grouped by their results on the diagnostic and post-tests, respectively. We used independent two-sample t-tests (equal variances not assumed) to determine if these groups means differ significantly.  We note that this analysis is moot for problem PT:5 as most students answered incorrectly.

\begin{center}
\begin{table}[h]
\caption{Comparison of mean final grade (MFG) of students by diagnostic question success}\label{MFGDiag}
    \begin{tabular}{c||c|c|c}
            Question & MFG (Correct/Justified)& MFG (Otherwise) &p-value \\
            \hline
         DT:1& 73.9&66.5 &0.072  \\
         DT:2& 81.2& 66.7& 0.019 \\
         DT:3& 75.4& 66.1&0.029  \\
         DT:4& 70.0&67.9 & 0.582 \\
         DT:5& * & * & * \\
    \end{tabular}
    \end{table}
\end{center}

\begin{center}
\begin{table}[h]
\caption{Comparison of mean final grade (MFG) of students by post-test question success}\label{MFGPost}
    \begin{tabular}{c||c|c|c}
            Question & MFG (Correct/Justified)& MFG (Otherwise) &p-value \\
            \hline
         PT:1& 78.0&64.2 & $<$0.001 \\
         PT:2& 70.0& 67.1& 0.471 \\
         PT:3& 82.3& 65.8& $<$0.001  \\
         
    \end{tabular}
    \end{table}
\end{center}
\vspace{0.3cm}

 To aid in answering the question of whether learning logic and counterexample affects student success, for each of the parts of DT:1-DT:4 on the Math 2903 diagnostic test we looked at the mean final grades for the students who did not answer the question completely correctly (not both correct and justified). We then compared the mean final grades of those answering the PT:2 on the post-test question correctly (CR: correct and justified) and those that did not. We again performed an independent two-sample t-test (equal variances not assumed) to determine if these groups means differ significantly. The results are in the table below.\\

\begin{center}
\begin{table}[h]
\caption{Comparison of mean final grade (MFG) of students with incorrect diagnostic test answers by post-test question success}\label{MFGComp}
\begin{tabular}{c||c|c|c}
            Incorrect & MFG (PT:2 CR)& MFG (Otherwise) &p-value \\
            \hline
         DT:1& 69.1&61.6  &0.104  \\
         DT:2& 67.2& 65.9& 0.76 \\
         DT:3& 68.7& 60.2&0.053  \\
         DT:4& 68.9&66.1 & 0.610 \\
         
    \end{tabular}
    \end{table}
\end{center}

\section{Conclusions}

We set out to study the effect of minimal logical training on student outcomes. Our data shows interesting trends associated to the effect of understanding the idea of counterexamples in logic.  Through our investigations we found that certain problems types can be used as checks for students less likely to find success in their mathematics course. We now examine these ideas specifically in relation to our data and then once again in the context of our perspective on learning mathematics in general.

From the data above, we can reasonably conclude two important things: first, that very few students have a strong idea of counterexample entering Math 2903, though many do have some intuition about whether a statement is true or false, without the skills to definitively show why. Indeed, from Table \ref{Wason} it can be seen that only 2 out of 60 students chose the correct answer, and the plurality of students checking both P true and Q true.  This is in line with the findings presented in \cite{IS}. Additionally, only 8 of 62 students answered DT:2 correctly. %\textcolor{blue}{I don't know what you mean here specifically.  Is this in reference to the Wason problem?  Maybe we could use the 8/62 who got DT:2? SE: I tried to make this more clear.}
  
From the midterm data, it is worth noting that these pre-service teachers can learn logical reasoning and counterexample. Given the data for the post-test, it is clear that students that do successfully learn the idea of logic and counterexamples have more success in the course overall than students that do not.  From Table \ref{MFGPost} we note that while the Wason question (PT:1) was worth approximately 1 per cent of the students' final grade, students that do answer the question correctly do, on average, 13.8 percentage points better in the course than those that do not. We emphasize that this question tests exactly one logical idea, without any computational or algebraic skills needed. Indeed, the Wason Selection Task seems to be a highly discriminant question. More research is needed to determine whether this logical knowledge is a driver of higher course grades, or whether stronger students simply learn these ideas more quickly than other students.

  Question PT2, where the student needed to supply a counterexample, was not very predictive. We believe this is the case since the question may have been too intuitive, though students who did not understand counterexamples in the beginning of the course and did indeed answer this question correctly seemed to have a significant, or near-significant, difference in mean grades as opposed to students who did not learn the idea of counterexample. We conjecture that it is important to understand the idea of counterexample at its core and not just intuitively for easy cases, as a central idea of mathematics is using logical and mathematical techniques to understand non-intuitive problems. Indeed, if we chose a slightly more challenging counterexample problem it is possible that we would have seen more significant differences in Table \ref{MFGComp}.
  
  Question PT3 was highly predictive. We are not surprised by this data;  the skills to write a correct proof mean that you understand the statement and the associated objects it deals with in generality, and have the computational and algebraic skills necessary to complete this proof.
  
  In our comparison of mean final grades of students that did not answer the first four diagnostic questions correctly, we find that, while the mean final grades of the students who did answer PT:2 correctly was higher for all diagnostic test questions, we find that at $\alpha = 0.05$ none of these means were significantly higher than the groups of students that answered PT:2 incorrectly. However, we would conjecture that we may have found a significant difference in mean final grade if PT:2 was a more difficult question. 
  
  Putting this all together, we believe that the data supports the current research that counterexample plays a very strong role in mathematical understanding in early undergraduate mathematics. Indeed, a student gaining understanding of the idea significantly improves student success in their respective course.

%We find that the Wason question is very highly discriminating in terms of final grades. We note that the 

%\begin{itemize}
%    \item Students who are already good at logic (especially proof) do better in the course.
%    \item Basically no one understands the Wason question at first
%    \item After instruction, the Wason question seems to be highly predictive of success in the course
%    \item The proof question also seems to be predictive, though there were many students who did well %without getting 6b correct.
%    \item Given the two Wason points above, it seems that students that learn Wason style logic have better %outcomes than students that do not.
%    \item Some questions were not predictive, especially 6a, since it may have been too easy of a question. 
%\end{itemize}

\bibliographystyle{plain}

\end{document}